\input amstex
\voffset=-3pc

\documentstyle{amsppt}
\voffset=-3pc

\def\bC{\Bbb C}
\def\bM{\Bbb M}
\def\bR{\Bbb R}

\def\sH{\Cal H}
\def\cK{\Cal K}
\def\cO{\Cal O}

\def\var{\varphi}

\def\sa{{\text{sa}}}

\def\asa{A_{\text{sa}}}

\def\pr{\text{pr}}

\def\b1{\bold 1}

\magnification=\magstep1
\parskip=6pt
\NoBlackBoxes
\topmatter
\title Convergence of Functions of Self--Adjoint Operators and Applications
\endtitle  
\rightheadtext{Convergence and Applications}
\author Lawrence G.~Brown
\endauthor
%\endtopmatter

\abstract{The main result (roughly) is that if $(H_i)$ converges weakly to $H$ and if also $f(H_i)$ converges weakly to $f(H)$, for a single strictly convex continuous function $f$, then $(H_i)$ must converge strongly to $H$.
One application is that if $f({\text{pr}}(H))={\text{pr}} (f(H))$, where ${\text{pr}}$ denotes compression to a closed subspace $M$, then $M$ must be invariant for $H$.
A consequence of this is the verification of a conjecture of Arveson, that Theorem 9.4 of [Arv] remains true in the infinite dimensional case.
And there are two applications to operator algebras.
If $h$ and $f(h)$ are both quasimultipliers, then $h$ must be a multiplier.
Also (still roughly stated), if $h$ and $f(h)$ are both in $p\asa p$, for a closed projection $p$, then $h$ must be strongly $q$--continuous on $p$.}
\endabstract
\endtopmatter

\noindent 2013 Mathematics Subject Classification:  Primary 47B15,47A60, 46L05

\noindent
Keywords and phrases:\ self--adjoint operator, weak convergence, strong convergence, strictly convex function, Korovkin type theorem, Kaplansky density theorem, quasimultiplier, $q$--continuous

\bigskip
{\bf1.\ Introduction.}

Although the main result of this paper involves only elementary operator theory, the original motivation was from a couple of technical operator algebraic questions.
In [B1, Proposition 2.59(a)] it was shown that if $f$ is a non--linear operator convex function on an interval I, if $h$ is a self--adjoint quasimultiplier of a $C^*$--algebra $A$ such that $\sigma(h)\subset I$, where $\sigma(h)$ is the spectrum, and if $f(h)$ is also a quasimultiplier, then $h$ is in fact a multiplier of $A$.
And in [B2, Theorem 4.14(i)] it was shown (for $f$ as above) that if $p$ is a closed projection in $A^{**}$, the bidual, and if $h$ and $f(h)$ are both in $p\asa p$, then $h$ (which is still assumed self--adjoint) must be strongly $q$--continuous on $p$ (provided also that either $p$ is compact or $0\in I$ and $f(0)=0$).
(In general $S_{\sa}=\{x\in S\colon x^*=x\}$.)
In both cases ad hoc methods could be used to generalize the result to certain non--operator convex functions $f$, and the motivating questions were to find for which $f$ the results are true.

It turns out that both results are true for an arbitrary strictly convex continuous function $f$.
Also the best approach is to start with an elementary result about convergence in the weak operator topology.
This result, roughly stated in the abstract, is Theorem 2.1 below.
The application to compressions is Corollary 2.7, and the proof of Arveson's conjecture is Corollary 2.8.
The operator algebraic results are Theorems 3.1 and 3.6.  (Theorem 3.6 is not strictly speaking an application of Theorem 2.1, but its proof is modeled on that of 2.1.)
Remark 2.6 and Proposition 2.10 discuss converses to Theorem 2.1, and Remark 2.6 also discusses extensions of the theorem and relations to the operator--theoretic version of the Kaplansky density theorem ([Kap, Thm. 2] for the forward direction and [Kad, Cor.~3.7] for the converse, see also [Ped, \S 2.3]).

There are other relationships involving the results of this paper, and the full story of these relationships does not seem clear.
K.~Davidson suggested that our results have the flavor of a Korovkin type theorem and also pointed out the relationship to [Arv].
The relationship of Arveson's results to Korovkin type theorems is discussed in [Arv, Remarks 1.4 and 1.8].
Also, Theorem 2.1, in the special case where the interval $I$ occurring there is compact, could be deduced from either Corollary 2.7 or Theorem 3.1, and Corollary 2.7 can be deduced from Theorem 3.6. 

We begin with some elementary results, none of which are probably new, but the only reference we know is that Corollary 1.2 can be deduced from [BL, Lemma 3.4].
Corollary 1.2 is used in \S2 and Corollary 1.3 in \S3.
The set of bounded linear operators on a Hilbert space $\sH$ is denoted by $B(\sH)$.

\proclaim{Lemma 1.1}Let $\sH=\sH_1\oplus \ldots\oplus\sH_n$ be a direct sum of Hilbert spaces, $v=t_1 u_1 \oplus\ldots\oplus t_u u_n$ a unit vector in $\sH$ where $t_i\geq 0$ and $\|u_i\|=1$, and $P$ a positive operator in $B(\sH)$.
Then $(Pv,v)\leq\sum_1^n (Pu_i,u_i)$.
\endproclaim 

\demo{Proof}$P$ can be represented by an $n\times n$ matrix $(P_{ij})$, where $P_{ij}\in B(\sH_j,\sH_i)$.
If $w=u_1\oplus\ldots\oplus u_n$ and $Q$ is the operator with matrix $((\delta_{ij}-t_i t_j) P_{ij})$, then the conclusion is just the statement that $(Qw,w)\geq 0$.
Since the matrix $(\delta_{ij}-t_i t_j)$ is positive, this follows from the fact that the Hadamard product of positive matrices is positive.
\enddemo

\proclaim{Corollary 1.2}Let $\sH=\sH_1\oplus\ldots\oplus\sH_n$ be a direct sum of Hilbert spaces, $P$ a positive operator in $B(\sH)$, and let $P_{ii}$ be the compression of $P$ to $\sH_i$.
Then $\|P\|\leq \sum_1^n \|P_{ii}\|$.
\endproclaim

\demo{Proof}Let $v$ be a unit vector in $\sH$, and use the notation of the lemma.
Then $(Pv,v)\leq\sum_1^n (Pu_i,u_i)=\sum_1^n (P_{ii} u_i,u_i)\leq\sum_1^n \|P_{ii}\|$.
\enddemo

If $A$ is a $C^*$--algebra, its bidual $A^{**}$ is a von Neumann algebra, called the enveloping von Neumann algebra of $A$.
Bounded linear functionals on $A$ are also regarded as weak$^*$--continuous linear functionals on $A^{**}$.
A state $\varphi$ on $A$ is said to be supported by a projection $p$ in $A^{**}$ if $\varphi(\b1-p)=0$ where $\b1$ is the identity of $A^{**}$.
For $q,r\in A^{**}$ and $\psi$ a bounded linear functional on $A$, $q\psi r$ denotes the functional $a\mapsto \psi(qar)$, which is again weak$^*$--continuous on $A^{**}$.
The set of positive elements of $A$ is denoted by $A_+$.

\proclaim{Corollary 1.3}Let $\varphi$ be a state on a $C^*$--algebra $A$, and $p=p_1+\ldots+p_n$, where $p_1,\ldots,p_n$ are mutually orthogonal projections in $A^{**}$.
If $\varphi$ is supported by $p$, then $\varphi\leq\sum_1^n (1/\varphi(p_i))p_i\varphi p_i$, where the $i$'th term is taken to be 0 if $\varphi(p_i)=0$.
\endproclaim

\demo{Proof}We may assume $\varphi(p_i) > 0,\ \forall i$.
If not, replace $p$ by $p'=\sum' p_i$, where the terms for which $\varphi(p_i)=0$ are omitted in the new sum, and note that $\varphi$ is supported by $p'$.

Now let $(\sH_\varphi,\pi,v)$ be the result of the GNS construction for $\varphi$, and let $\tilde\pi$ be the canonical extension of $\pi$ to $A^{**}$.
We apply the lemma with $\sH=\tilde\pi (p)\sH_\varphi,\sH_i=\tilde\pi(p_i)\sH_\varphi$, and $P$ the compression of $\pi(a)$ to $\sH$ for $a$ in $A_+$, noting that the state vector $v$  is indeed in $\sH$.
The conclusion is that $\varphi(a)=(Pv,v)\leq\sum_1^n (\pi(a) u_i,u_i)=\sum_1^n (1/\varphi (p_i))(p_i\varphi p_i)(a)$.
\enddemo

{\bf 2.\ The main result and basic applications.}

\proclaim{Theorem 2.1}Let $f$ be a strictly convex continuous function on an interval $I$, let $\sH$ be a Hilbert space, and let $(H_i)$ be a net in $B(\sH)_{\sa}$ such that $\sigma(H_i)\subset I$, $\forall i$.
If $(H_i)$ converges weakly to a bounded operator $H$ with $\sigma(H)\subset I$, and if also $(f(H_i))$ converges weakly to $f(H)$, then $(\varphi(H_i))$ converges strongly to $\varphi(H)$ for every bounded continuous function $\varphi$ on $I$.
In particular if the net $(H_i)$ is bounded, then $(H_i)$ converges strongly to $H$.
\endproclaim

{\bf{Note.}}  The applications described in the introduction depend only on the special case where the interval $I$ is compact.

\medskip\noindent
{\bf Beginning of the proof}.
Since it is sufficient to prove that every subnet of $(H_i)$ has a further subnet for which the conclusion is true, we may replace the given net with a subnet.
For a Borel set $E$, let $P_i(E)=\chi_E (H_i)$ and $P(E)=\chi_E(H)$, the spectral projections.
Choose the subnet so that for each interval $E$, the net $(P_i(E))$ converges weakly to a positive contraction $Q(E)$.
Note that $Q(\cdot)$ is finitely additive.
We are going to prove a close relationship between $Q(\cdot)$ and $P(\cdot)$ via several lemmas.

Let $a_0$ be the minimum point in $\sigma(H)$ and $b_0$ the maximum point.
If $a_0$ is the left endpoint of $I$, let $a=a_0$.
Otherwise choose a point $a$ in $I$ such that $a<a_0$.
Similarly, $b=b_0$ if $b_0$ is the right endpoint of $I$ and $b\in I$, $b>b_0$ otherwise.
Thus in all cases $[a,b]$ is a compact subinterval of $I$ and $\sigma(H)\subset [a,b]$.

The first lemma is needed to deal with the endpoints, and it also illustrates one of the main ideas of the proof.

\proclaim{Lemma 2.2}If $a<x<y<b$, then 
$$
\| P( [a,x]) Q([y,\infty)) P([a,x])\| \leq {x-a\over y-x}.
$$
\endproclaim

\demo{Proof}Let $g$ be the function obtained by subtracting a linear function from $f$ so that $g(a)=g(x)=0$.
Note that $(g(H_i))$ converges weakly to $g(H)$.
Let $-\gamma=g(c)$ be the minimum value of $g$.
Because of strict convexity, $\gamma>0$ and $a<c<x$.
Also $g$ is strictly increasing to the right of $x$ and in particular $\beta=g(y)>0$.
For each $i$ $\beta P_i ([y,\infty))\leq g(H_i)+\gamma \b1$, where $\b1$ is the identity operator on $\sH$.
Thus $\beta Q([y,\infty))\leq g(H)+\gamma \b1$, and
$$
\align
\beta P([a,x]) Q([y,\infty)) P([a,x])&\leq P([a,x]) g(H) P([a,x])+\gamma P([a,x])\\
&\leq\gamma P([a,x]).
\endalign
$$
Hence
$$
\| P([a,x]) Q([y,\infty)) P([a,x])\|\leq \gamma/\beta.\tag1
$$
Now the one--sided derivatives $g'(z\pm)$ exist and are finite, $\forall z\in (a,b)$.
By strict convexity, $\gamma/(x-a) < \gamma/(x-c) < g' (x-)$, and $\beta/(y-x) > g'(x+) \geq g' (x-)$.
Thus $\gamma/\beta\leq {x-a\over y-x}$.
\enddemo

The next lemma is a standard real analysis fact stated in our notation.
We need $x_0>a$ because of the possibility that $f'(a+)=-\infty$.

\proclaim{Lemma 2.3}If $a< x_0 < x < b$ and if $\epsilon >0$, then there is a partition $x_0 < x_1 < \ldots <x_n=x$ such that $f'(x_j-) -f' (x_{j-1}+) <\epsilon$ for $j=1,\ldots,n$.
\endproclaim

\demo{Proof}There are only finitely many points $y$ in $(x_0,x)$ such that $f'(y+)-f'(y-)\geq \epsilon$.
Let these points be $y_1,\ldots,y_m$ from left to right.
For each interval $J=[x_0,y_1]$, $[y_{k-1},y_k]$, or $[y_m,x]$, there is $\delta > 0$ such that $z,w\in J$ and $|z-w|<\delta$ implies $f'(w-) -f' (z+) < \epsilon$.
This is proved in the same way as the uniform continuity of continuous functions on $J$.
So we can obtain the desired partition as a suitable refinement of $x_0 < y_1 < \ldots <y_n < x$.
\enddemo

\proclaim{Lemma 2.4}If $a<x<y<b$, then $P([a,x]) Q([y,\infty))=0$.
\endproclaim

\demo{Proof}Let $\eta=f' ({x+y\over 2}\ +)-f' (x-) > 0$.
Choose $\epsilon>0$, and choose $x_0$ in $(a,x)$ such that $x_0-a<\epsilon$.
Let $x_0<x_1<\ldots <x_n=x$ be a partition as in Lemma 2.3.

For each $j=1,\ldots,n$ let $g_j$ be the function obtained by subtracting a linear function from $f$ so that $g_j(x_j-1)=g_j (x_j)=0$.
Let $-\gamma_j=g(c_j)$ be the minimum value of $g_j$, and let $\beta_j=g_j (y)$.
As in the proof of (1) above, we see that 
$$
\|P([x_{j-1},x_j]) Q([y,\infty)) P([x_{j-1},x_j])\|\leq \gamma_j/\beta_j.
$$
Since $g'_j (x_j-) -g'_j (x_{j-1}+)=f' (x_j-) -f' (x_{j-1} +) < \epsilon$, and since $g'_j (x_j -) > 0$ and $g'_j (x_{j-1} +) < 0$, then $g'_j (x_j-) < \epsilon$.
Thus $\gamma_j/(x_j-c_j) < g'_j (x_j-) < \epsilon$, and $\gamma_j < \epsilon (x_j-c_j) < \epsilon (x_j-x_{j-1})$.
Also, since $g'_j\ ({x+y\over 2}\ +)-g'_j(x-)$$=f'\ ({x+y\over 2}\ +)-f'(x-)=\eta$, and since $g'_j(x-)\geq g'_j (x_j-)>0$, then $g'_j\ ({x+y\over 2}\ +) > \eta$.
Thus $\beta_j > g_j\ ({x+y\over 2})+\eta\ ({y-x\over 2}) > g_j (x_j)+\eta\ ({y-x\over 2})=\eta ({y-x\over 2})$.
So $\|P([x_{j-1},x_j]) Q([y,\infty)) P([x_{j-1},x_j])\|\leq {\epsilon\over\eta}\ {2\over y-x}\ (x_j-x_{j-1})$.
Also, by Lemma 2.2, $\|P([a,x_0]) Q[y,\infty))P([a,x_0])\|$
$$
< {\epsilon\over y-x_0}\leq {\epsilon\over y-x}.
$$

Now Corollary 1.2 implies that
$$
\gathered
\|P([a,x]) Q([y,\infty)) P([a,x])\|\leq \|P([a,x_0]) Q([y,\infty)) P([a,x_0])\| +\\
\sum_1^n \|P([x_{j-1},x_j]) Q([y,\infty)) P([x_{j-1},x_j\|\\
\leq\ {\epsilon\over y-x}+{\epsilon\over \eta}\ {2\over y-x}\ (x-a).
\endgathered
$$
Since $\epsilon$ is arbitrary, we conclude that $P([a,x]) Q([y,\infty)) P([a,x])=0$, and since $Q([y,\infty))\geq 0$ this implies $P([a,x]) Q([y,\infty))=0$.
\enddemo

\proclaim{Lemma 2.5}If $a<y<b$, then
$$
P((-\infty,y))\leq Q((-\infty,y))\leq Q((-\infty,y])\leq P((-\infty,y]).
$$
Also $P([a,b])=Q([a,b])=\b1$.
\endproclaim

\demo{Proof}Since $P((-\infty,y))=P([a,y))=\lim\limits_{x\to y-} P([a,x])$, Lemma 2.4 implies $P((-\infty,y)) Q([y,\infty))=0$.
This implies $Q([y,\infty))\leq P([y,\infty))$ and $Q((-\infty,y))=\b1-Q([y,\infty))\geq \b1 -P([y,\infty))=P((-\infty,y))$.
Using left--right symmetry, we can also prove $Q((-\infty,y])\leq P((-\infty,y])$ and $Q((y,\infty))\geq P((y,\infty))$.
If $b$ is the right endpoint of $I$, then $Q((b,\infty))=0$, since $P_i((b,\infty))=0$, $\forall i$.
Otherwise, choose $x$ and $y$ so that $b_0<x<y<b$ and apply Lemma 2.4.
Since $P([a,x])=1$, we conclude that $Q([y,\infty))=0$, whence $Q((b,\infty))=0$.
A symmetrical proof shows that $Q((-\infty,a))=0$.
\enddemo

\noindent
{\bf End of the proof}.
It is enough to show that $(\varphi((H_i))$ converges weakly to $\varphi (H)$ for each bounded continuous function $\varphi$ on $I$, since then we also get the same conclusion for $|\varphi|^2$, and the facts that $(\varphi(Hi))$ converges weakly to $\varphi(H)$ and $(\varphi(H_i)^2)$ converges weakly to $|\varphi (H)|^2$ imply that $(\varphi(H_i))$ converges strongly to $\varphi(H)$.
Then fix $\varphi$ and a unit vector $v$ in $\sH$.
Let $\varepsilon >0$ and choose $\delta$ such that $x,y\in [a,b]$ and $|x-y|<\delta$ imply $|\varphi(x)-\varphi(y)|<\epsilon$.
Note that the monotone function $\alpha(x)=(P((-\infty,y))v,v)$ has only countably many discontinuities in $(a,b)$, and at each continuity point, $(P((-\infty,x)) v,v)=(P((-\infty,x])v,v)$.
Choose a partition $a=x_0<x_1<\ldots < x_n=b$ of mesh $<\delta$ such that $x_j$ is a continuity point of $\alpha$ for $j=1,\ldots,n-1$.
Let $E_j=[x_{j-1},x_j)$ for $1\leq j < n$ and $E_n=[x_{n-1},x_n]$.
Finally, let $T=\sum\limits_1^n \varphi (x_{j-1}) P(E_j)$ and let $T_i=$
$$
\varphi (H_i) P_i ((-\infty,a)) +\sum_1^n\varphi (x_{j-1}) P_i (E_j)+\varphi (H_i) P_i ((b,\infty)).
$$
Then $\|T-\varphi (H)\|<\epsilon$ and $\|T_i-\varphi (H_i)\|<\epsilon$.

Now let $s$ be any cluster point of the net $((\varphi (H_i) v,v))$.
Note that the bounded nets $(P_i ((-\infty,a)))$ and $(P_i ((b,\infty)))$ converge strongly to 0, since they are positive and converge weakly to 0.
Therefore $((T_i v,v))$ converges to $(Tv,v)$.
Thus $|s-(Tv,v)|\leq \epsilon$, and clearly $| (\varphi (H) v,v)- (Tv,v)|\leq \epsilon$.
Since $\epsilon$ is arbitrary, we conclude $s=(\varphi (H) v,v)$, and hence $(\varphi (H_i))$ converges weakly to $\varphi (H)$.

\example{Remarks 2.6}(i)\ If the net is not bounded, then $(H_i)$ may not converge strongly to $H$.

Let $1<\alpha<\beta\leq 2$, and let $(K_i)$ be a net in $B(\sH)_{\sa}$ such that $(K_i)$ converges strongly to $K$ but $(|K_i|^\alpha)$ does not converge strongly to $|K|^\alpha$.  Then $(|K_i|^{\alpha\over 2})$ converges strongly to $|K|^{\alpha\over 2}$ and $(|K_i|^{\beta\over 2})$ converges strongly to $|K|^{\beta\over 2}$.
Therefore $(|K_i|^\alpha)$ converges weakly to $|K|^\alpha$ and $(|K_i|^\beta)$ converges weakly to $|K|^\beta$, and we get a counterexample by taking $H_i=|K_i|^\alpha$ and $f(x)=|x|^{\beta/\alpha}$.  Both [Kap] and [Kad] are used to justify the above.

(ii)\ The Kaplansky density theorem can be deduced from Theorem 2.1.
If $(H_i)$ is a net in $B(\sH)_{\sa}$ which converges strongly to $H$, then $(H_i^2)$ converges weakly to $H^2$.
Thus the hypotheses of 2.1 are met with $I=\bR$ and $f(x)=x^2$, and the conclusion that $(\varphi(H_i))$ converges strongly to $\varphi(H)$ for $\varphi$ bounded and continuous is sufficient for the Kaplansky density theorem.

(iii)\ We sketch an argument for how to fill the gap between what was proved in the theorem and what was disproved in (i).
Assume the hypotheses of the theorem and let $g$ be a function, obtained by subtracting a linear function from $f$, such that $g(a)>0$ and $g'(a-)<0$ if $a$ is not the left endpoint of $I$ and $g(b)>0$ and $g'(b+)>0$ if $b$ is not the right endpoint of $I$.  (If $a_0<b_0$, the most natural choice is to take $g(a_0)=g(b_0)=0$; and if $a_0=b_0$ and is in the interior of $I$, the most natural choice is to make the $x$-axis a line of support for g at $a_0$.) 
The arguments at the end of the proof, with $g$ in place of $\varphi$, show that$(g(H_i)P_i([a,b]))$ converges weakly to $g(H)$.  
Since also $(g(H_i))$ converges weakly to $g(H)$, we conclude that $(g(H_i) (P_i(-\infty,a)+P_i (b,\infty))$ converges weakly to 0.
So if $\varphi$ is a continuous function on $I$ such that $|\varphi|/g(x)$ is bounded on $I\backslash [a,b]$, the arguments in the proof show that $(\varphi(H_i))$ converges weakly to $\varphi(H)$.  In particular 
if $x^2/g(x)$ is bounded on $I\backslash [a,b]$, the weak convergence of $(H_i^2)$ to $H^2$ implies that $(H_i)$ converges strongly to $H$.
These conditions  can be restated, without mentioning $g$, $a$, or $b$, as follows:\ \ For $\var (H_i)$ to converge weakly to $\var (H)$, it is sufficient that $\var(x) / {\text{max}}(1, |x|, f(x))$ be bounded on $I$.  And for $(H_i)$ to converge strongly to $H$, it is sufficient that if $I$ is unbounded on the right, then $f(x)>0$ for $x$ sufficiently large and $x^2=\cO (f(x))$ as $x\to\infty$; and if $I$ is unbounded on the left, then $f(x)>0$ for $x$ sufficiently small and $x^2=\cO (f(x))$ as $x\to -\infty$.

The argument in [Kad] can be adapted to show that the above results are sharp in the following sense:\ If for all nets and a particular $f$ and $\var$  the hypotheses of the theorem imply that $(\var(H_i))$ converges weakly to $\var(H)$, then the first condition above must be satisfied; and if for all nets and a particular $f$  the hypotheses imply that $(H_i)$ converges strongly to $H$, then the second  condition must be satisfied. 

(iv)  It is possible to weaken the hypotheses of the theorem:  Instead of assuming that $(f(H_i))$ converges weakly to $f(H)$, assume only that $(f(H)v,v)\geq {\text {lim sup}}(f(H_i)v,v)$ for each $v$ in $\sH$.  The same proof works.  We do not know whether this strengthening of the theorem is valuable, but note that it is irrelevant if $f$ is operator convex and $I$ is compact, as in the result cited from [B1].  In this case it is automatic that $(f(H)v,v) \leq{\text{lim inf}}(f(H_i)v,v)$.  This follows from [B1, Proposition 2.34].  (See [B3, Theorem 3.1 and remark following] for an elementary treatment.)
\endexample

\proclaim{Corollary 2.7}Let $H\in B(\sH)_{\sa}$, let $M$ be a closed subspace of the Hilbert space $\sH$, and let $f$ be a strictly convex continuous function on an interval that includes $\sigma(H)$.
If $f({\text{pr}} (H))=\text{ pr}(f(H))$, where $\pr$ denotes compression to $M$, then $M$ is invariant for $H$.
\endproclaim

\demo{Proof}Let $x_0$ be a point in $\sigma(H)$.
We may assume $M^\perp$ is infinite dimensional.
If not, replace $\sH$ with $\sH\oplus\ell^2$, $H$ with $H\oplus x_0 \b1$, and $M$ with $M\oplus 0$.
Let $S$ be a unilateral shift on $M^\perp$ (of uncountable multiplicity if $M^\perp$ is non--separable), and let $V=\b1_M\oplus S$, regarding $\sH$ as $M\oplus M^\perp$.

Then $(V^n H V^{*n} +x_0 (\b1- V^n V^{*n}))$ converges weakly to $\pr (H) \oplus x_0 \b1_{M^\perp}$.
Since 
$$
f(V^n H V^{*n} +x_0 (\b1-V^n V^{*n}))=V^n f(H) V^{*n} +f(x_0) (\b1-V^n V^{*n}),
$$
then $(f(V^n HV^{*n}+x_0 (\b1-V^n V^{*n}))$ converges weakly to $\pr (f(H))+f(x_0) \b1_{M^\perp}=f(\pr (H)\oplus x_0 \b1_{M^\perp})$.
So the theorem implies that $(V^n H V^{*n} +x_0 (\b1-V^n V^{*n}))$ converges strongly to $\pr (H) \oplus x_0 \b1_{M^\perp}$.
If $u$ is a vector in $M$ and $Hu= u' \oplus w$, this implies that $\|S^n w\|\to 0$ and hence $w=0$.
\enddemo

Corollary 2.7 leads to the elimination of the finite dimensionality hypothesis in a theorem of Arveson.  The special case where $f$ is operator convex had previously been proved by Petz.

\proclaim{Corollary 2.8}(cf.~[Arv, Theorem 9.4], [Pet]) Let $f$ be a strictly convex continuous function on a compact interval $[a,b]$, $\varphi\colon B(\sH)\to B(\cK)$ a unital completely positive map where $\sH$ and $\cK$ are Hilbert spaces, and let $H\in B(\sH)_{\sa}$ with $\sigma(H)\subset [a,b]$.
If $\varphi (f(H))=f(\varphi (H))$, then the restriction of $\varphi$ to the algebra of polynomials in $H$ is multiplicative.
\endproclaim

\demo{Proof}By the Stinespring theorem there are a Hilbert space $\tilde\cK \supset \cK$ and a unital $*$--representation $\pi\colon B(\sH)\to B(\tilde\cK)$ such that $\varphi (T)=\pr (\pi (T))$, $\forall T\in B(\sH)$.
Since $\pi (f(H))=f(\pi (H))$, Corollary 2.7 implies that $\cK$ is invariant for $\pi (H)$.
\enddemo

\proclaim{Corollary 2.9}If $f$ is a strictly convex function on a compact interval $[a,b]$ and $\epsilon >0$, then $\exists \delta >0$ such that:

For any Hilbert space $\sH$, any closed supspace $M$, and any $H$ in $B(\sH)_{\sa}$ with $\sigma(H)\subset [a,b]$, 
$\| f(\pr (H))-\pr (f(H))\|<\delta\Rightarrow \| PH-HP\| < \epsilon$, where $P$ is the projection with range $M$ and $\pr$ denotes compression to $M$.
\endproclaim

\demo{Proof}If not, then for each $n=1,2,\ldots$, there are $\sH_n,M_n$, and $H_n$ such that $\sigma(H_n)\subset [a,b]$, $\|f(\pr (H_n))-\pr (f (H_n))\| < {1\over n}$ and $\|P_n H_n-H_n P_n\|\geq \epsilon$.
Let $\sH=\bigoplus_1^\infty \sH_n$, $M=\bigoplus_1^\infty M_n$, a closed subspace of $\sH$, $H=\bigoplus_1^{\infty} H_n$, an element of $B(\sH)_{\sa}$ with $\sigma(H)\subset [a,b]$, and $P=\bigoplus_1^{\infty} P_n$, the projection with range $M$.
Further, let $A$ be the $C^*$--subalgebra of $ B(\sH)$ consisting of operators $\bigoplus_1^\infty T_n$ with  $\{\|T_n\|\}$ bounded, I the closed two--sided ideal of $A$ consisting of operators $\bigoplus_1^\infty T_n$ with $\| (T_n)\|\to 0$, and $\pi\colon A\to A/I$ the quotient map.
If $h=\pi(H)$ and $p=\pi (P)$, then $pf(h) p=pf (php+a (\b1-p))p$ but $ph\neq hp$.
Since $A/I$ can be faithfully represented on a Hilbert space, this contradicts Corollary 2.7.
\enddemo

Finally, we prove a converse to Theorem 2.1.

\proclaim{Proposition 2.10}If $f$ is a continuous real--valued function on a compact interval $[a,b]$ which is neither strictly convex nor strictly concave, then there are an operator $H$ and a sequence $(H_n)$ in $B(\ell^2)_{\sa}$ with $\sigma(H)$, $\sigma(H_n)\subset [a,b]$, such that $(H_n)$ converges weakly to $H$ and $(f(H_n))$ converges weakly to $f(H)$, but $(H_n)$ does not converge strongly to $H$.
\endproclaim

\demo{Proof}By the proof of [Arv, Proposition 9.2] there are $x,y\in [a,b]$ and $t\in (0,1)$ such that $x<y$ and $f(tx+(1-t)y)=tf(x)+(1-t) f(y)$.
Let $H_0$ be the $2\times 2$ matrix 
$$
x \pmatrix t&\sqrt{t(1-t)}\\ \sqrt{t(1-t)}&1-t\endpmatrix + y\pmatrix 1-t&-\sqrt{t(1-t)}\\ -\sqrt{t(1-t)}&t\endpmatrix.
$$
Then $\sigma (H_0)=\{x,y\}\subset [a,b]$, and $f(\pr (H_0))=\pr f(H_0)$ but $M$ is not invariant for $H_0$, where $M$ is the one--dimensional subspace of $\bC^2$ consisting of vectors of the form $\pmatrix *\\ 0\endpmatrix$.
The proof of Corollary 2.7 produces the required counterexample.
\enddemo

{\bf {3.\ Operator algebraic applications.}}

If $A$ is a $C^*$--algebra, let $M(A)=\{ t\in A^{**}\colon tA+At\subset A\}$, the algebra of multipliers of $A$, and let $QM(A)=\{t\in A^{**}\colon AtA\subset A\}$, the set of quasimultipliers of $A$.
Also $S(A)$ denotes the state space of $A$ and $Q(A)$ denotes $\{\varphi\in A^*\colon \varphi\geq 0$ and $\|\varphi\|\leq 1\}$, the quasi-state space of $A$.

\proclaim{Theorem 3.1}Let $f$ be a continuous strictly convex function on a compact interval $[a,b],\ A$ a $C^*$--algebra, and $h$ an element of $A^{**}_{\sa}$ such that $\sigma(h)\subset [a,b]$.
If both $h$ and $f(h)$ are in $QM (A)$, then $h\in M (A)$.
\endproclaim

\demo{Proof}We will show that $h^2\in QM (A)$.
The result then follows from [AP, Proposition 4.4].
By [AP, Theorem 4.1], it is sufficient to show $\varphi_i (h^2)\to\varphi (h^2)$ whenever $\varphi_i$ and $\varphi$ are in $S(A)$ and $\varphi_i\to\varphi$ weak$^*$.
We may assume that the GNS representations for $\varphi$ and the $\varphi_i$'s are all in Hilbert spaces of the same dimension.
If not, replace $A$ by $A\otimes\cK(\sH)$, where $\cK(\sH)$ is the algebra of compact operators on a Hilbert space $\sH$ of sufficiently large dimension, so that $A^{**}$ is replaced by the $W^*$--tensor product $A^{**}\Bar{\otimes} B(\sH)$.
Then replace $h$ by $h\otimes \b1$, $\varphi$ by $\varphi\otimes\psi$, and $\varphi_i$ by $\varphi_i\otimes\psi$, for a pure state $\psi$ on $\cK(\sH)$.

Now, reverting to the original notation, and passing to a subnet, which is permissible, we can realize all the GNS representations by maps $\pi_i,\pi\colon A\to B(\sH)$ with the same unit vector $v$ as the state vector such that $\pi_i (a)\to \pi(a)$ strongly, $\forall a\in A$, (cf.~[Dix, Section 3.5]).
Let $\tilde\pi_i$ and $\tilde\pi$ denote the canonical extensions to $A^{**}$.
If $t\in QM (A)$ and $a,b\in A$, the facts that $\pi_i (a)^* \tilde\pi_i (t) \pi_i (b)\to \pi (a)^*\tilde\pi (t)\pi(b)$ strongly, $\pi_i(a)\to\pi(a)$ strongly, and $\pi_i(b)\to \pi(b)$ strongly imply that $(\tilde\pi_i (t)\pi(b) u,\pi (a) w)\to (\tilde\pi (t)\pi(b) u$, $\pi(a) w), \forall u,w\in\sH$.
Since $\pi$ is non--degenerate, this implies $\tilde\pi_i(t)\to\pi(t)$ weakly.
So Theorem 2.1 now applies with $H_i=\tilde\pi_i (h)$.
Thus $\tilde\pi_i (h)\to\tilde\pi (h)$ strongly and hence $\tilde\pi_i (h^2)=\tilde\pi_i (h)^2\to\tilde\pi (h)^2 =\tilde\pi (h^2)$ weakly (also strongly).
Therefore $\varphi_i (h^2)=(\tilde\pi_i (h^2) v,v)\to (\tilde\pi (h^2) v,v)=\varphi (h^2)$.
\enddemo

\example{Remark} As in Remark 2.6(iv), we can weaken the hypotheses:  Instead of assuming $f(h)$ is in $QM(A)$, assume only that $f(h)$ is weakly upper semicontinuous.  This strengthening of the theorem is irrelevant if $f$ is operator convex, since then, by [B1, Proposition 2.34], $f(h)$ is automatically weakly lower semicontinuous.
\endexample

A converse to Theorem 3.1 follows from [B1, Example 2.66] but is not stated clearly there.

\proclaim{Proposition 3.2}Let $f$ be a continuous real--valued function on the compact interval $[a,b]$ which is neither strictly convex nor strictly concave, and let $A$ be the $C^*$--algebra of norm convergent sequences in $\cK(\ell^2)$.
Then there is $h$ in $A_{\sa}^{**}$ such that $\sigma(h)\subset [a,b]$, $h$ and $f(h)$ are in $QM (A)$, and $h$ is not in $M(A)$.
\endproclaim

\demo{Proof}It is easy to see that $A^{**}$ can be identified with the algebra of bounded collections $t=\{ T_n\colon 1\leq n\leq \infty, T_n\in B(\ell^2)\}$.
Then $t$ is in $QM (A)$ if and only if $t_n\to t_\infty$ weakly and $t\in M(A)$ if and only if $t_n\to t$ strongly.
So the result follows from Proposition 2.10.
\enddemo

For a projection $p$ in $A^{**}$, $F(p)$ denotes $\{\varphi\in Q(A)\colon \varphi(\b1-p)=0\}$, the face of $Q(A)$ supported by $p$.
Then $p$ is called {\it closed} [Ak1] if $F(p)$ is weak$^*$ closed and {\it compact} ([Ak2]) if $F(p)\cap S(A)$ is weak$^*$ closed.
If $p$ is closed and $h\in p A_{\sa}^{**} p$, then $h$ is called {\it strongly lower semicontinuous on} $p$ ([B2]) if the map $\varphi  \mapsto\varphi (h)$ is lower semicontinuous (lsc) on $F(p)$, and $h$ is called {\it strongly $q$--continuous on} $p$ ([B1], cf.~also [APT]) if $\chi_F (h)$ is closed whenever $F$ is a closed subset of $\bR$ and compact if also $0\not\in F$.
Here $\chi_F (h)$ denotes the spectral projection computed in $p A^{**} p$.
If $A$ is unital the qualifier ``strongly'' is unnecessary and every closed projection is compact.
It was shown in [B1, Theorem 3.43] that $h$ is strongly $q$--continuous on $p$ if and only if $h=pa$ for some $a$ in $A_{\sa}$ such that $ap=pa$.

\proclaim{Lemma 3.3}Let $A$ be a unital $C^*$--algebra, $p$ a closed projection in $A^{**}$, and $f$ a strictly convex continuous function on a compact interval $[a,b]$.
If $h\in p A_{\sa}p$ such that $\sigma(h) \subset [a,b]$ and $f(h)$ is (strongly) lower semicontinuous on $p$, then $h$ is (strongly) $q$--continuous on $p$.
\endproclaim

The proof follows the next two lemmas, in which $p(E)$ denotes $\chi_E (h)$, the spectral projection computed within $p A^{**} p$ and the notation of Lemma 3.3 is assumed.

\proclaim{Lemma 3.4}Let $a<x<y<b$.
If the $\varphi_i$'s are states of $A$ supported by $p([a,x])$ and $\varphi_i\to\varphi$ weak$^*$, then $\varphi (p([y,b]))\leq {x-a\over y-x}$.
\endproclaim

\demo{Proof}Let $g,c,\gamma$, and $\beta$ be as in the proof of Lemma 2.2
Then $\beta p ([y,b])\leq g(h) +\gamma p$ and $g(h)$ is lsc on $p$.
Therefore
$$
\beta \varphi (p ([y,b]))\leq \varphi (g(h)) +\gamma\leq \liminf \varphi_i (g (h))+\gamma\leq \gamma.
$$
Hence
$$
\varphi(p ([y,b]))\leq \gamma/\beta.\tag2
$$
As in the proof of Lemma 2.2, $\gamma/\beta\leq {x-a\over y-x}$.
\enddemo

\proclaim{Lemma 3.5}Let $a< x<y<b$.
If the $\varphi_i$'s are states of $A$ supported by $p([a,x])$ and $\varphi_i\to\varphi$ weak$^*$, then $\varphi (p([y,b]))=0$.
\endproclaim

\demo{Proof}Let $\eta$ be as in the proof of Lemma 2.4, let $\epsilon > 0$, and choose $x_0,\ldots,x_n$ as in the proof of 2.4.
By Corollary 1.3, for each $i$ we can write $\varphi_i\leq\sum_0^n\psi_{ij}$, where $\psi_{io}$ is a state of $A$ supported by $p([a,x_0])$ and $\psi_{ij}, j\geq 1$, is a state of $A$  supported by $p([x_{j-1},x_j])$.
Passing to a subnet, we may assume each net $(\psi_{ij})$ converges weak$^*$ to a state $\varphi_j$.
Clearly $\varphi\leq\sum_0^n\varphi_j$.
By Lemma 3.4,
$$
\varphi_0 (p([y,b]))\leq {x_0-a\over y-x_0} < {\epsilon\over y-x_0} \leq {\epsilon\over y-x}.
$$
\enddemo

As in the proof of Lemma 3.4, for $j\geq 1$, $\varphi_j (p([y,b]))\leq \gamma_j/\beta_j$, in the notation of Lemma 2.4, and as in the proof of 2.4, 
$\gamma_j/\beta_j\leq {\epsilon\over\eta}\ {2\over y-x}\ (x_{j}-x_{j-1})$.
Therefore $\varphi(p([y,b]))\leq {\epsilon\over y-x} +\sum^n_1\ {\epsilon\over\eta}\ {2\over y-x}\ (x_j-x_{j-1})\leq {\epsilon\over y-x}+{\epsilon\over\eta} {2\over y-x} (x-a)$.
So since $\epsilon$ is arbitrary, $\varphi (p([y,b]))=0$.

\demo{Proof of Lemma 3.3}Let $a<x<b$, and let the $\varphi_i$'s be states of $A$ supported by $p([a,x])$ such that $\varphi_i\to\varphi$ weak$^*$.
Choose a monotone sequence $(y_n)$ such that $y_n\in (x,b)$ and $y_n\to x$.
Since $p([y_n,b])\to p((x,b))$ in the weak* topology of $A^{**}$ and since $\varphi(p[y_n,b])=0$, $\forall n$, then $\varphi(p(x,b])=0$.
In other words $\varphi$ is supported by $p([a,x])$, whence $p([a,x])$ is closed.
A symmetrical proof shows that $p([x,b])$ is closed.
Akemann showed in [Ak1] that the set of closed projections is closed under arbitrary lattice meets and also that $q_1 \vee q_2$ is closed if $q_1$ and $q_2$ are and if the angle between $q_1$ and $q_2$ is positive.
The latter applies in particular if $q_1q_2=q_2q_1$.
So we can conclude that $p(E)$ is closed for all closed sets $E$.
\enddemo

\proclaim{Theorem 3.6}Let $A$ be a $C^*$--algebra, $p$ a closed projection in $A^{**}$, and $f$ a strictly convex continuous function on a compact interval $[a,b]$.
Assume also that either $p$ is compact or $0\in [a,b]$ and $f(0)=0$.
If $h\in p A_{\sa} p$ such that $\sigma(h)\subset [a,b]$ and $f(h)$ is strongly lower semicontinuous on $p$, then $h$ is strongly $q$--continuous on $p$.
Here $\sigma(h)$ and $f(h)$ are computed within $p A^{**}p$.
\endproclaim

\demo{Proof}We apply Lemma 3.3 to $\tilde A$, the result of adjoining an identity to $A$, identifying $\tilde A^{**}$ with $A^{**}\oplus\bC$.
In the compact case, we use $p\oplus 0$ in place of $p$ and $h\oplus 0$ in place of $h$.
Since $f(h\oplus 0)=f(h)\oplus 0$, the conclusion is immediate.
\enddemo

In the non--compact case, we use $\tilde p=p\oplus 1$ in place of $p$ and $h\oplus 0$ in place of $h$.
It is still true that $f(h\oplus 0)=f(h)\oplus 0$, and we conclude that $h\oplus 0$ is (strongly) $q$--continuous on $\tilde p$.
If $E$ is a closed set not containing 0, then $\chi_E (h\oplus 0)=\chi_E (h) \oplus 0$; and the fact that this is closed in $\tilde A^{**}$ implies that $\chi_E(h)$ is compact in $A^{**}$.
If $E$ is a closed set that contains 0, then $\chi_E (h\oplus 0)=\chi_E (h)\oplus 1$; and the fact that this is closed in $\tilde A^{**}$ implies that $\chi_E(h)$ is closed in $A^{**}$.
Therefore $h$ is strongly $q$--continuous on $p$.

\example{Remarks 3.7}(i)\ (cf. Remark 2.6(iv) and the remark following Theorem 3.1) It follows from the theorem that $f(h)$ is in $p A_{\sa} p$.
We could have proved a slightly weaker theorem by assuming $f(h)$ in $p A_{\sa}p$ instead of $f(h)$ strongly lsc, and we do not know whether the extra strength of the actual theorem is worthwhile.
In [B2], this issue did not arise because $f$ was assumed operator convex.
Theorem 4.3 of [B2] implies in this case that $f(h)$ is automatically strongly usc on $p$.

(ii)\ Remark 4.15 of [B2] gives some discussion, which will not be repeated here, of the hypothesis that $0\in [a,b]$ and $f(0)=0$.
We mention only that if $p$ is not compact, it is impossible to have $f(0)<0$, given the other hypotheses.
\endexample

\proclaim{Proposition 3.8}Let $f$ be a continuous real--valued function on a compact interval $[a,b]$.
If $f$ is neither strictly convex nor strictly concave, then there are a unital $C^*$--algebra $A$, a closed projection $p$ in $A^{**}$, and an $h$ in $p A_{\sa}p$ such that $\sigma(h)\subset [a,b]$, $f(h)\in pAp$ and $h$ is not $q$--continuous on $p$.
\endproclaim

\demo{Proof}Let $A$ be the algebra of convergent sequences in $\bM_2$, the algebra of $2\times 2$ matrices.
Then $A^{**}$ can be identified with the set of bounded collections $t=\{ T_n\colon 1\leq n\leq\infty,\ T_n\in\bM_2\}$.
Let the projection $p$ be given by $p_n=\pmatrix 1&0\\ 0&0\endpmatrix$ for $n<\infty$ and $p_\infty=\pmatrix 1&0\\ 0&1\endpmatrix$.
Then $p$ is closed.
As in the proof of proposition 2.10, there is a self--adjoint matrix $H=\pmatrix \alpha&\beta\\ \overline\beta&\gamma\endpmatrix$ such that 
$\sigma(H)=\{x,y\}\subset [a,b]$, $x<\alpha<y$, and $f(H)=\pmatrix f(\alpha)&*\\ *&*\endpmatrix$.
Let $h$ be given by $H_n=\pmatrix \alpha&0\\ 0&0\endpmatrix$ for $n<\infty$ and $H_\infty=H$.
Then $h$ and $f(h)$ are in $p A_{\sa}p$.
If $E$ is the closed set $\{\alpha\}$, then $\chi_E(h)$ is not closed.
\enddemo

\bye